\documentclass[12pt,leqno]{amsart}

\usepackage[utf8]{inputenc}   % LaTeX, comprends les accents !
\usepackage[T1]{fontenc}      % Police contenant les caract??res français
\usepackage{geometry}         % Définir les marges
\usepackage[francais]{babel}  % Placez ici une liste de langues, la
\usepackage{amssymb}
\usepackage{enumerate}
\voffset-1cm \markleft{\rm Poisson} 
\newtheorem{theorem}{Theorem}
\newtheorem{definition}[theorem]{Definition}

\newtheorem{lemma}[theorem]{Lemma}

\newtheorem{proposition}[theorem]{Proposition}

\newtheorem{defi}{D\'{e}finition}
\newtheorem{propo}{Proposition}

\newcommand{\K}{\mathbb K}

\newcommand{\g}{\frak{g}}

\newcommand{\N}{\mathbb N}
\newcommand{\R}{\mathbb R}

\newcommand{\p}{\mathcal{P}}
\newcommand{\KS}{{\mathbb K}[\Sigma_3]}
\newcommand{\ST}{\Sigma_3}

\newcommand{\ra}{\rightarrow}

\newcommand{\pf}{\noindent{\it Proof. }}

\newcommand{\wa}{weakly associative }
\newcommand{\im}{{\rm Im }}

\usepackage{graphicx}

\title{Weakly associative algebras, Poisson algebras and quantization}
\author{Elisabeth Remm}
\date{2 sivan 5780}
\address{ Universit\'e de Haute Alsace, Laboratoire IRIMAS-Math\'ematiques, 6 rue des FR\`eres Lumi\`ere. F.68093 Mulhouse}
\email{elisabeth.remm@uha.fr}

\begin{document}

\maketitle

\noindent{INTRODUCTION} 

Given a commutative associative algebra $(A,\mu)$ on a field $\K$ of characteristic $0$, any associative formal deformation $(A[[t]],\mu_t)$ with $\mu_t=\mu+ \sum_k t^k\varphi _k$  of this algebra determines a Poisson algebra $(A,\mu,\psi)$ and $(A[[t]],\mu_t)$ is a quantization by deformation of this Poisson algebra. Recall that the formal deformation $\mu_t$ is associative if the following identity is formally satisfied
$$\mu_t(X,\mu_t(Y,Z))-\mu_t(\mu_t(X,Y),Z)= 0.$$
Developing this identity, we obtain
\begin{enumerate}
\item Order $0$: $\mu$ is associative, this is given in the hypothesis,
\item Order $1$: 
 $$\mu(X,\varphi_1(Y,Z))-\varphi_1(X,\mu(Y,Z))-\mu(\varphi_1(X,Y),Z)-\varphi_1(\mu(X,Y),Z)= \delta^2_{H,\mu}\varphi_1 (X,Y,Z)=0$$ where $\delta^2_{H,\mu}$ is the coboundary operator associated with the Hochschild cohomology of $(A,\mu)$,
\item Order $2$ : $\varphi_1(X,\varphi_1(Y,Z))-\varphi_1(X,\varphi_1(Y,Z))+ \delta^2_{H,\mu}\varphi_2 (X,Y,Z)=0$
\end{enumerate}
One of the fundamental problems is to know, given a Poisson algebra $(A,\mu,\psi)$, if there is a formal deformation of the associative commutative algebra $(A,\mu)$ which is a deformation by quantification of the Poisson algebra  $(A,\mu,\psi)$. Let us take, for example, the algebra $\mathcal{C}^\infty(M)$ of $\mathcal{C}^\infty$-functions on a differential manifold $M$, with the associative commutative product $\mu(f,g)=fg.$ Let us assume that we have also a Poisson structure $(A,\mu,\psi)$ on this associative commutative algebra, that is a Lie bracket $\psi$ which satisfies 
$$\psi(f,gh)=g\psi(f,h)+\psi(f,h)g.$$
In this case, Maxim Kontsevich has proved (\cite{Ko}) that there exists a formal associative deformation $\mu_t$ of $\mu$ which is a deformation quantization of $(A,\mu,\psi)$. 

Now, still within the framework of the quantization and  deformation of associative commutative algebras, we can consider the problem of whether Poisson algebras can be obtained by more general deformations, that is to say in deforming associative algebra in a class of properly chosen class of nonassociative algebras. In this work we show that such a class exists, the algebras of this class are called weakly associative. More generally, we study deformations  of weakly associative algebras and we define a generalization of Poisson algebras by considering the nonassociative case. We  study the algebraic structure of this algebra class. 
Recall also that an algebra $(A,\mu)$ is a nonassociative algebra  if it is a vector space over K and is equipped with a bilinear binary multiplication operation $\mu: A \times A \rightarrow  A$ which may or may not be associative.
%\tableofcontents

\section{Weakly associative algebras}

Let $\K$ be a field of characteristic $0$ and $(A,\mu)$ be a $\K$-nonassociative algebra. To make the writing of this multiplication $\mu$ more natural , we will write in this paragraph and the following, $\mu(X,Y)=X \ast Y$ and the algebra $(A,\ast)$. Thus the associator of the algebra $(A,\ast)$ is the trilinear map  $\mathcal{A}_\ast$ defining by
$$\mathcal{A}_\ast (X,Y,Z)=X \ast (Y \ast Z)- (X \ast Y )\ast Z$$
for any $X,Y,Z \in A.$

\subsection{Definition and Examples}

\begin{defi}
The $\K$-algebra $(A,\ast)$ is called weakly associative if its associator $\mathcal{A}_\ast$ satisfies
$$\mathcal{A}_\ast (X,Y,Z)+\mathcal{A}_\ast (Y,Z,X)-\mathcal{A}_\ast (Y,X,Z)=0$$
for any $X,Y,Z \in A.$
\end{defi}
To simplify, we shall denote by $\mathcal{WA}_\ast (X,Y,Z)$ the trilinear map
$$\mathcal{WA}_\ast (X,Y,Z)=\mathcal{A}_\ast (X,Y,Z)+\mathcal{A}_\ast (Y,Z,X)-\mathcal{A}_\ast (Y,X,Z)$$
and a $\K$-algebra $(A,\ast)$ is weakly associative  if and only if
$$\mathcal{WA}_\ast (X,Y,Z)=0$$
for any $X,Y,Z \in A.$
\medskip

\noindent{\bf Examples.}
\begin{enumerate}
  \item Any associative  algebra is weakly associative.
  \item Any abelian (i.e. commutative) algebra is weakly associative.
  \item Any Lie algebra is weakly associative. In fact, let $(\g, \{,\})$ be a $\K$-Lie algebra.  The associator of the Lie bracket satisfies
  $$\mathcal{A}_{\{,\}} (X,Y,Z)=\{X,\{Y,Z\}\}-\{\{X,Y\},Z\}$$
  and from the Jacobi identity
  $$\mathcal{A}_{\{,\}}  (X,Y,Z)=\{\{Z,X\},Y\}.$$
  Then
  $$\mathcal{WA}_{\{,\}}  (X,Y,Z)=\{\{(Z,X\},Y\}+\{\{X,Y\},Z\}+\{\{Y,Z\},X\}$$
  and, from the Jacobi identity
   $$\mathcal{WA}_{\{,\}}  (X,Y,Z)=0$$
 \item Let $E$ be a set and $\mathcal{F}(E,\R)$ be the vector space of applications $f: E \rightarrow \R$. Consider $a,b \in \R$ with $ a \neq b$ and   the noncommutative multiplication on $E$ :
 $$f \star g= af +bg.$$
 Then $\mathcal{A}_\star (f,g,h)=(a^2-a)f+(b-b^2)g.$ If we assume that $a,b \neq 0$ then this operation is not associative. But
 $$\mathcal{WA}_\star (f,g,h)=(a^2-b^2-a +b)f.$$
 Then it is not commutative as soon as $b \neq a$ and weakly associative as soon as $b=1-a$.  However, we have no distributivity.
\item In a same way, for any $a \neq 0$,  the binary operation $f \ast g =f^ag^a$ is not associative and weakly associative. 
 
  \item Let $\widehat{\K[X,Y]}$ be the graded vector space of polynomials with two noncommutative indeterminates $X$ and $Y$ and let $I$ the ideal generated by 
   $$
   \left\{
   \begin{array}{l}
   (XX)X-X(XX),\ \ (YY)Y-Y(YY), \\
   (XY)X-X(YX), \ \ (YX)Y-Y(XY), \\
   (YX)X-X(XY)+(XX)Y-Y(XX), \ \ (XY)Y-Y(YX)+(YY)X-X(YY).
   \end{array}
   \right.
   $$
  The algebra $\widehat{\K[X,Y]}/I$ is weakly associative. 

\item An example in dimension $2$.
Let us consider on  $\K^2=\K\{e_1,e_2\}$ the multiplication corresponding to
$$
\left\{
\begin{array}{l}
 e_1e_1=\frac{a}{2}e_1,   \\
e_1e_2=\frac{a+2}{4}e_2, \\   `
e_2e_1=\frac{a-2}{4}e_2, \\   
e_2e_2=0.
\end{array}
\right.
$$
Since $\mathcal{A}(e_1,e_1,e_2)=\frac{a^2-4}{16}e_2,\mathcal{A}(e_2,e_1,e_1)=-\frac{a^2-4}{16}e_2$ and $\mathcal{A}(e_i,e_j,e_k)=0$ for all the other cases, this product is weakly associative. It is associative if and only if $a=2$ or $a=-2$. For example,
$$(x_1,y_1)(x_2,y_2)=(3x_1x_2,2x_1y_2+y_1x_2)$$
is weakly associative. 

\subsection{Free weakly associative algebras}
Let us denote by $(\mathcal{F}_{WA}(X),\ast)$ the free weakly associative algebra with $1$-generators. It is a graded algebra
$$\mathcal{F}_{WA}(X)= \oplus_{k \geq 0} \mathcal{F}_{WA}(X)(k)$$
where $\mathcal{F}_{WA}(X)(k)$ consists of the vector space of elements of degree $k$. We put $$\mathcal{F}_{WA}(X)(0)=\K.$$ Let us determine each  of these components. For this, we use the relations
\begin{equation}\label{aaa}
\mathcal{WA}_\ast (a,a,a)=\mathcal{A}_\ast(a,a,a),\mathcal{WA}_\ast (a,a,b)=\mathcal{WA}_\ast (a,b,a)=\mathcal{A}_\ast(a,b,a).
\end{equation}
\begin{enumerate}
  \item Degree $2$: We put $X \ast X= X^2$. Thus $\mathcal{F}_{WA}(X)(2)=\K\{X^2\}.$
  \item Degree $3$: We have $\mathcal{WA}_\ast (X,X,X)=X\ast X^2-X^2 \ast X=0.$ We put $X^3=X^2\ast X=X \ast X^2$ and $\mathcal{F}_{WA}(X)(3)=\K\{X^3\}.$
  \item Degree $4$: $$\left\{
  \begin{array}{l}\mathcal{WA}_\ast (X,X,X^2)=\mathcal{WA}_\ast (X,X^2,X)=X \ast X^3-X^3 \ast X \\
  \mathcal{WA}_\ast (X^2,X,X)=0.
  \end{array}
  \right.
  $$Then we put
  $$X^3 \ast X=X \ast X^3=X^4_1, \ \ X^2 \ast X^2 = X^4_2,$$
  and $\mathcal{F}_{WA}(X)(4)=\K\{X^4_1,X^4_2\}.$
  \item Degree $5$: $\mathcal{WA}_\ast (X,X,X^3)=X \ast X^4_1-X^4_1 \ast X$, $\mathcal{WA}_\ast (X^3,X,X)=X^3 \ast X^2-X^2 \ast X^3.$  Likewise
  $\mathcal{WA}_\ast (X^2,X^2,X)=\mathcal{WA}_\ast (X^2,X,X^2)=X^2 \ast X^3-X^3 \ast X^2$ and $\mathcal{WA}_\ast (X,X^2,X^2)=-X^4_2 \ast X+X \ast X^4_2.$ We put
  $$
  \left\{
  \begin{array}{l}
     X^3 \ast X^2=X^2 \ast X^3=X^5_1,    \\
      X^4_1 \ast X=X \ast X^4_1 =X^5_2, \\
     X^4_2 \ast X=X \ast X^4_2 =X^5_3 \\
\end{array}
\right.
$$
and  $\mathcal{F}_{WA}(X)(5)=\K\{X^5_1,X^5_2,X^5_3\}.$
\end{enumerate}
Assume by induction, that for a given degree $n$, we have $u\ast v=v\ast u$ for any $u \in \mathcal{F}_{WA}(X)(p)$ and  $v \in \mathcal{F}_{WA}(X)(q), p+q \leq n$. The relations (\ref{aaa}) imply that $(a\ast b)\ast a=a\ast (b\ast a)$ and $u\ast v=v\ast u$ for any $u \in \mathcal{F}_{WA}(X)(p)$ and  $v \in \mathcal{F}_{WA}(X)(q)$, with $ p+q \leq n+1.$ Then
$$\mathcal{F}_{WA}(X)(n)= \oplus_{p+q=n} \widetilde{\mathcal{F}_{WA}(X)(p)\otimes \mathcal{F}_{WA}(X)(q)}$$
where $\widetilde{}$ denotes the symmetric quotient.  In particular
$$d_{2p+1}=\dim \mathcal{F}_{WA}(X)(2p+1)= \sum _{k=1}^{p}d_kd_{2p+1-k}$$
and
$$d_{2p}=\dim \mathcal{F}_{WA}(X)(2p)= \sum _{k=1}^{p-1}d_kd_{2p-k}+ \frac{d_p(d_p+1)}{2}$$
with $d_1=d_2=1.$
\end{enumerate}

\subsection{Some properties of weakly associative algebras}

\begin{proposition}\label{der}
Let $(A,\ast)$ be a nonassociative $\K$-algebra. Then for any $X \in A$, the endomorphism
$L_X-R_X$ defined by
$$(L_X-R_X)(Y)=X \ast Y - Y \ast X$$
is a derivation of $(A,\ast)$ if and only if $(A,\ast)$ is weakly associative.
\end{proposition}
\pf Recall that a derivation $f$ of $(A,\ast)$ is a linear endomorphism which satisfies
$$f(X) \ast Y + X \ast f(Y)=f(X \ast Y)$$
for any $X,Y \in A$. 
Then
$$
\begin{array}{l}
  (L_X-R_X)(Y) \ast Z+Y \ast (L_X-R_X)(Z)-(L_X-R_X)(Y \ast Z)   \\ 
   \ \ \ = (X \ast Y)\ast Z+Y \ast (X \ast Z)-X \ast( Y \ast Z) -(Y \ast X) \ast Z - Y \ast (Z \ast X)+ (Y \ast Z)\ast X\\
   \ \ \  = -\mathcal{A}_\ast (X,Y,Z)-\mathcal{A}_\ast (Y,Z,X)+\mathcal{A}_\ast(Y,X,Z)\\
  \ \ \ = -\mathcal{WA}_\ast (X,Y,Z).
\end{array}
$$
and $L_X-R_X$ is a derivation of $(A,\ast)$ if and only if $\mathcal{WA}_\ast=0.$

For example, if $(A,\ast)$ is a Lie algebra, then $L_X-R_X=2adX$ and we find again that $ad X$ is a derivation.

\medskip
\begin{proposition}
Any weakly associative algebra is a Lie admissible algebra and also a flexible algebra.  
\end{proposition}
\pf Recall that $(A,\ast)$ is Lie admissible if the multiplication $\psi(X,Y)=X \ast Y- Y \ast X$ is a Lie bracket.  This is equivalent to the identity (\cite{GRLieadm})
$$\mathcal{A}_\ast (X,Y,Z)-\mathcal{A}_\ast (Y,X,Z)-\mathcal{A}_\ast (Z,Y,X)-
\mathcal{A}_\ast (X,Z,Y)+\mathcal{A}_\ast (Y,Z,X)+\mathcal{A}_\ast (Z,X,Y)=0.$$
This identity is equivalent to
$$\mathcal{WA}_\ast (X,Y,Z)-\mathcal{WA}_\ast (X,Z,Y)=0$$
implying that any weakly associative algebra is Lie admissible.
The algebra $(A,\ast)$ is flexible if 
$$\mathcal{A}_\ast (X,Y,X)=0.$$
This is equivalent to
$$\mathcal{A}_\ast (X,Y,Z)+\mathcal{A}_\ast (Z,Y,X)=0.$$
But
$$\mathcal{WA}_\ast (Z,X,Y)+\mathcal{WA}_\ast (Y,X,Z)=\mathcal{A}_\ast (X,Y,Z)+\mathcal{A}_\ast (Z,Y,X)$$
then the weakly associativity implies the flexibility. The converse is not true. Note however that free weakly associative algebra with one generator coincides with the free flexible algebra with one generator. But this is no longer true for $2$ or more generators.

\medskip

\noindent{\bf Remark.} Let $LIE$ be the category of Lie algebras. The objects are the Lie algebras, and the morphisms the Lie algebra homomorphisms. We can also consider the category $WASS$ of weakly associative algebra. We define the functor $T$ from $LIE$ to $WASS$ which consists to associate to a Lie algebra $(A,\mu)$ the weakly associative  algebra $(A,\mu)$. If $f:(A,\mu) \ra (B,\mu')$ is a morphism of Lie algebra, that is $f(\mu(X,Y))=\mu'(f(X),f(Y)),$ we have the same identity for $T(f)$ considering $\mu$ as a weakly associative multiplication. Then this functor $T$ is well defined. Conversely, if $(A,\ast)$ is a weakly associative algebra, since it is Lie admissible, the product $\mu_\ast(X,Y)=(X \ast Y - Y \ast X)$ is a Lie bracket on $A$ and $(A, \mu)$ is a Lie algebra. The application $T_2$ which associate to a weakly associative algebra $(A,\ast)$ the Lie algebra $(A,\mu_\ast)$ is a functor between the categories $WASS$ and $LIE$. But we have not an equivalence of categories because $T_1 \circ T_2 (A,\ast)=(A,\mu_\ast)$. Let us remark also that $LIE$, $ASS$ and $WASS$ are subcategories of the category $ALG$ whose objects are the $\K$-algebras whose multiplication is given by a bilinear map satisfying a quadratic relation. We  have  noticed that $LIE$ and $ASS$ are subcategories of $WASS$. In this context, $WASS$ occurs as the smallest subcategory of $ALG$ containing $LIE$ and $ASS$.

\subsection{$\Sigma_3$-associative algebras, $v$-algebras}
In \cite{GRnonass} we introduce a class of nonassociative algebras  whose associator satisfies some symmetric relations. We shall see that weakly associative algebras also belong to this class.  Let  $\Sigma_3=\{Id,\tau_{12},\tau_{13},\tau_{23},c,c^2\}$ be the symmetric group of degree $3$ where $\tau_{ij}$ is the transposition between $i$ and $j$ and $c$ the cycle $(231)$. We denote  by $\K[\Sigma_3]$ the algebra group of $\Sigma_3$. It is provided with an associative algebra structure and with a $\Sigma_3$-module structure. The (left) action of $\ST$ on  $\KS$ is given by
$$(\sigma \in \ST, v=\sum a_i\sigma_i \in \KS) \ra \sum a_i \sigma_i \circ \sigma.$$ For any $v \in \KS$, $\mathcal{O}(v)$ is the corresponding orbit and $\K[\mathcal{O}(v)]$ the $\K$-linear space of $\KS$ generated by $\mathcal{O}(v)$. It is also a $\ST$-submodule. For example, let us consider the vector of $\KS$,  $v_{WA}=Id+c_1-\tau_{12}$. Its orbit contains the vectors
$$
 \left\{
 \begin{array}{ll}
 v_{WA} \circ Id& =  Id+c-\tau_{12} \\
 v_{WA} \circ \tau_{12}    &  = \tau_{12} +\tau_{23}-Id  \\
     v_{WA}\circ \tau_{13}    &  = \tau_{13} +\tau_{12}-c  \\
    v_{WA} \circ \tau_{23}    &  = \tau_{23} +\tau_{13}-c^2  \\
   v_{WA} \circ c   &  = c+c^2-\tau_{13}  \\
    v_{WA} \circ c^2   &  = c^2+Id-\tau_{23}  \\ 
\end{array}
  \right.
  $$
 The $\ST$-invariant subspace $\K[ \mathcal{O}( v_{WA})]$ of $\KS$ is of dimension $4$ admitting $\{ Id+c-\tau_{12},\tau_{12} +\tau_{23}-Id,\tau_{13} +\tau_{1}-c, \tau_{23} +\tau_{13}-c^2\}$ as a linear basis.

\noindent{\bf Notation} Let $A$ be a $\K$-vector space. We consider, for any $\sigma \in \ST$ the linear map $\Phi^A_\sigma$ or $\Phi_\sigma$  if there is no ambiguity on A:
$$\Phi^A_\sigma : A ^{\otimes^3} \ra A ^{\otimes^3}$$
given by
$\Phi^A_\sigma (X_1 \otimes X_2 \otimes X_3)=X_{\sigma(1)} \otimes X_{\sigma(2)} \otimes X_{\sigma(3)}$. For any $v=\sum a_i\sigma_i \in \KS$, we put  $\Phi^A_v=\sum a_i\Phi^A_{\sigma_i}.$ 

\begin{definition} \cite{GRnonass} An algebra $(A,\ast)$ is called $\KS-associative$ if there is $v \in \KS$ such that the quadratic relation of definition of $A$ is $\mathcal{A}_\ast \circ  \Phi^A_v=0.$ We shall say also that $(A,\ast)$ is a $v$-algebra.
\end{definition}
\begin{proposition}
Any  weakly associative algebras  satisfies $$\mathcal{A}_\ast \circ  \Phi^A_{v_{WA}}=0$$ with $v_{WA}=Id+c-\tau_{12}$. 
It is a $\KS$-associative algebra.
\end{proposition}
Of course the relation of definition $\mathcal{A}_\ast \circ  \Phi^A_v=0$ is not unique. For example a weakly associative algebra is a $v_{WA}$-algebra. But it is also a $(\tau_{12}-Id+\tau_{23})$-algebra. We shall say that two quadratic relations $\mathcal{A}_\ast \circ  \Phi^A_v=0$ and $\mathcal{A}_\ast \circ  \Phi^A_{v'}=0$ are equivalent if $ \K[\mathcal{O}(v)]=\K[\mathcal{O}(v')]$. In \cite {GRnonass}, we have classified all the relations of type $\mathcal{A}_\ast \circ  \Phi^A_v=0$. In particular  weakly associative algebra coincides with the type $IV(1)$ (\cite{GRnonass}, Theorem{\bf 1}) which writes
$$\mathcal{C}_\ast(\alpha)(x,y,z)=2\mathcal{A}_\ast (x,y,z)+(1+\alpha)\mathcal{A}_\ast (y,x,z)+\mathcal{A}_\ast (z,y,x)+\mathcal{A}_\ast (y,z,x)+(1-\alpha)\mathcal{A}_\ast (z,x,y) =0$$
with $\alpha =-1/2.$\medskip

\noindent{\bf Remark.} Assume that the multiplication $\ast$ is skew-symmetric and satisfies $C_{\ast}(\alpha)(x,y,z)=0.$  This identity  is equivalent to
$$ (1-2\alpha)(x \ast y)\ast z +2 (y \ast z)\ast x -2\alpha (z \ast x)\ast y =0.$$
We find again that $\ast$ is a Lie bracket if and only if $\alpha=-1/2$. In all the other cases this multiplication is of Lie type in the sense of \cite{Maka}.

\section{Nonassociative Poisson algebras}

\subsection{Definition}

\begin{definition}
Let  $(P,\bullet,\{,\})$ be  a triple where $P$ is a $\K$-vector space and $\bullet, \{,\}$  two multiplications on $P$. We say that $(P,\bullet,\{,\})$ is a nonassociative Poisson if
  \begin{enumerate}
  \item $(P,\bullet)$ is a  nonassociative commutative algebra,
  \item $(P,\{,\})$ is a Lie algebra,
	\item the Leibniz identity between $\{,\}$ and $\bullet$ is satisfied, that is
$$\{X\bullet Y,Z\}=X\bullet\{Y,Z\}+\{X,Z\}\bullet Y$$
fo all  $X,Y,Z \in P.$
\end{enumerate}
\end{definition}
\noindent{\bf Remark.}  As the terminology suggests, nonassociative Poisson algebras generalize Poisson algebras by relaxing the
associative condition on the underlying commutative algebra.

\subsection{Polarization of $\K$-algebras} The polarization technique, introduced in \cite{MRPoisson}, consists in representing a given one-operation $\K$-algebra $(A, \ast)$ without particular symmetry as an algebra with two operations, one commutative
and the other anticommutative. Explicitly, in the following, we will decompose the $\K$-algebra $(A, \ast)$ using:
\begin{enumerate}
  \item $X\bullet Y= X\ast Y + Y \ast X$ its symmetric part,
  \item $\{X,Y\}= X\ast Y - Y \ast X$ its skewsymmetric part,
\end{enumerate}
for $X,Y \in A$. 
The triplet $(A, \bullet, \{., .\})$  will be referred to as the polarization of $(A, \ast)$.
The polarization technique of \cite{MRPoisson} allows to make a link between weakly associative algebras and the notion of
nonassociative Poisson algebra.

\begin{theorem}
Let  $(P,\bullet,\{,\})$ be a nonassociative Poisson algebra. Let us consider on $P$ the third multiplication $$X \ast Y= X\bullet Y + \{X,Y\}.$$
Then the algebra $(P, \ast)$ is weakly associative.

%Conversely, if $(P,\ast)$ is a weakly associative algebra, then its polarized $(P, \bullet, \{.,.\}$ is a non-associative Poisson algebra.
\end{theorem}
\pf  We have to prove that $\mathcal{A}_\ast \circ \Phi^A_{v_{WA}}=0$ with $v_{WA}= Id+c_1-\tau_{12}.$  Since the 
multiplication  $\ast$ is Lie admissible,  we have $\mathcal{A}_\ast \circ \Phi^A_W =0$ with $W=Id-\tau_{12} - \tau_{13} -\tau_{23}+c_1+c_2$. The Leibniz identity writes  $\mathcal{A}_\ast \circ \Phi^A_w=0$ with
$w=Id+\tau_{12} + \tau_{13} -\tau_{23}-c_1+c_2.$ Then $W+w \in \K[\mathcal{O}(v_{WA})]$ this implies  $\mathcal{A}_\ast \circ \Phi^P_{v_{WA}}=0$ and $\ast$ is weakly associative.

\begin{theorem}
Let  $(P,\ast)$ be a weakly associative algebra. Then its polarization  $P(\bullet,\{,\})$ is a nonassociative Poisson algebra. We shall say also that $(P,\ast)$ is an admissible nonassociative Poisson algebra. 
\end{theorem}
\pf Since a weakly associative algebra is Lie admissible, the bracket   $\{X,Y\}=X \ast Y - Y \ast X$ is a Lie bracket. The multiplication $X \bullet Y=X \ast Y + Y \ast X$ is commutative. We have to show the Leibniz identity which is also written $\mathcal{A}_\ast \circ \Phi^P_w =0. $ But  this is a consequence of $w \in \K[\mathcal{O}(v_{WA})]$. 

\begin{lemma}
 Let  $( A,\ast)$ a Lie admissible algebra and $(A,\bullet, \{.,.\})$ its polarization. Then the following relations are equivalent
 \begin{enumerate}
  \item (Leibniz($\{,\},\cdot$): $\{X\ast Y,Z\}-X \ast \{Y,Z\}-\{X,Z\} \ast Y=0,$
  \item  (Leibniz($\{,\},\bullet$): $\{X\bullet Y,Z\}-X\bullet\{Y,Z\}-\{X,Z\}\bullet Y=0,$
\end{enumerate}
for any $X,Y,Z \in A.$
\end{lemma}
\pf Assume that
$$\{X\ast Y,Z\}-X\ast \{Y,Z\}-\{X,Z\}\ast Y=0$$
Since $2X\ast Y=\{X,Y\}+X\bullet Y$, we have
$$\{X\bullet Y,Z\}-X\bullet\{Y,Z\}-\{X,Z\}\bullet Y +\{\{X,Y\},Z\}-\{X,\{Y,Z\}\}-\{\{X,Z\},Y\}=0.$$
But $(A,\ast)$ is Lie admissible. Then   $\{.,.\}$ is a Lie bracket and
$$\{\{X,Y\},Z\}+\{\{Y,Z\},X\}+\{\{Z,X\},Y\}=0.$$
Then
$$\{X\bullet Y,Z\}-X\bullet\{Y,Z\}-\{X,Z\}\bullet Y.$$
The converse is similar.

\medskip

 Let $(A,\ast)$ be a weakly associative algebra and $(A,\bullet,\{.,.\})$ its polarization. We have seen that  any weakly associative algebra is Lie admissible. As consequence, the polarized bracket $\{.,.\}$ is a Lie bracket.  But we have no similar property for the commutative multiplication $\bullet$.  Recall that, if $(A,\ast)$ is an associative algebra, then the algebra $(A,\bullet)$ is a Jordan algebra, that is an abelian algebra such that
$$(X \bullet Y) \bullet X^2=X \bullet (Y \bullet X^2)$$
for any $X,Y \in P$.  In the weakly associative case, we have
$$\begin{array}{lll}
\mathcal{A}_\bullet(X,Y,Z) &=&(X\ast Y + Y \ast X) \ast Z + Z \ast (X \ast Y + Y \ast X) -X \ast (Y \ast Z + Z \ast Y)\\
&&-(Y \ast Z + Z \ast Y) \ast X \\
& =&\mathcal{A}_\ast(X,Y,Z)-\mathcal{A}_\ast(Z,Y,X)-(Y \ast Z)\ast X-X \ast(Z \ast Y)+Z\ast(X \ast Y)\\
&&+(Y\ast X)\ast Z.
\end{array}$$
This implies 
$$\mathcal{A}_\bullet(X,Y,Z)+\mathcal{A}_\bullet(Z,Y,X)=0$$
and $(A,\bullet)$ is a flexible algebra. Remember that this property is always true for commutative algebra. 
When $Z=X^2$, the previous equation becomes
$$\begin{array}{lll}
\mathcal{A}_\bullet(X,Y,X^2) &
 =&\mathcal{A}_\ast(X,Y,X^2)-\mathcal{A}_\ast(X^2,Y,X)-(Y \ast X^2)\ast X-X \ast(X^2 \ast Y)\\
&&+X^2\ast(X \ast Y)+(Y\ast X)\ast X^2.
\end{array}$$
Since $\mathcal{A}_\ast(X,X,X)=X \ast X^2 -X^2 \ast X=0$, then 
$$\mathcal{A}_\bullet(X,Y,X^2) =\mathcal{A}_\ast(X,Y,X^2)-\mathcal{A}_\ast(X^2,Y,X)=\mathcal{A}_\ast \circ \Phi_{Id-\tau_{13}} (X,Y,X^2).$$
But  $Id+\tau_{13}=(Id+\tau_{13}) \circ v_{WA}.$Thus $Id+\tau_{13} \in \mathcal{O}(v_{WA}).$ Since $\tau_{13}=(Id+\tau_{13})v_{WA}-Id$, we have                                                                  
$$\mathcal{A}_\bullet(X,Y,X^2)=2\mathcal{A}_\ast(X,Y,X^2).$$
We deduce 
\begin{proposition}
Let $(A,\ast)$ be a weakly associative algebra and $(A,\bullet,\{.,.\})$ its polarization. Then
 $(A,\bullet)$ is a commutative Jordan algebra if and only if $(A,\ast)$ is a Jordan algebra.
\end{proposition}

\section{On the cohomology of weakly associative algebras}
\subsection{The Lichnerowicz cohomology of a nonassociative Poisson algebra} A Poisson algebra is defined by two multiplications, a commutative associative algebra, a Lie bracket and these multiplications are linked by the Leibniz identity. The notion of deformation of Poisson algebra is connected directly to the notion of deformation quantization. Classically, one deforms the Lie bracket on a Poisson algebra considering that the associative product is preserved. We will recall in the next section the notion of formal deformation and we shall see that these deformations are parametrized by some cohomological complex. The complex associated with the Poisson bracket deformation is the Lichnerowicz cohomology. In this section we generalize the notion of this complex to the weakly associative case.

  Let $k$ be a non null integer. A $k$-cochain of the Lichnerowicz complex of a nonassociative Poisson algebra $(A, \bullet, \{.,.\},)$ is a skew-symmetric $k$-linear map
$$\Phi_k: A ^k \rightarrow A$$
which is a derivation on each arguments. Let us denote by  $\mathcal{C}^k(A)$ the vector space of $k$-cochains  and $\mathcal{C}^0(A)= A.$ The coboundary operator $\delta^k$ is the linear map
$$\delta^k: \mathcal{C}^k(A) \rightarrow \mathcal{C}^{k+1}(A)$$
defined by
$$
\begin{array}{ll}
\delta_L^k (\Phi_k)(X_0,X_1,\cdots,X_k)= &\sum_{i \geq 0} (-1)^i \{X_i,\Phi_k(X_0,\cdots,\widehat{X_i},X_{i+1},\cdots,X_k)\} \\
& \\
&+ \sum_{0 \leq i<j \leq k}(-1)^{i+j}\Phi_k(\{X_i,X_j\},X_0,\cdots,\widehat{X_i},\cdots,\widehat{X_j},\cdots,X_k).
\end{array}$$
The symbol $\widehat{X}$ means that the term is forgotten  For example
$$\delta_L^0 (X)(X_0)=\{X_0,X\}$$
$$\delta_L^1(\Phi_1)(X_0,X_1)=\{X_0,\Phi_1(X_1)\}-\{X_1,\Phi_1(X_0)\}-\Phi_1(\{X_0,X_1\})$$
and
$$
\begin{array}{ll}
\delta_L^2(\Phi_2)(X_0,X_1,X_2)=&\{X_0,\Phi_2(X_1,X_2)\}-\{X_1,\Phi_2(X_0,X_2)\}+\{X_2,\Phi_2(X_0,X_1)\} \\
&-\Phi_2(\{X_0,X_1\},X_2)+\Phi_2(\{X_0,X_2\},X_1)-\Phi_2(\{X_1,X_2\},X_0).$$
\end{array}
$$
Since
$$\delta_L^{k+1} \circ \delta^k_L$$
we have a cohomological complex.

\subsection{The weakly associative cohomology}

In this section we look how we can generalize the Hochschild cohomology concerning associative algebras to the weakly associative algebras. Let us begin to recall the expression of the Hochschild  operator.  
 Let $\phi$ be a $k$-linear map on the vector space $A$. The expression of this operator is
  $$
  \begin{array}{ll}
  \delta_H^k\phi(X_1,\cdots,X_{k+1})=&X_1\phi(X_2,\cdots,X_{k+1})+\sum_{i=1}^{k+1}(-1)^{i}\phi(X_1,\cdots,X_iX_{i+1},\cdots,X_{k+1})\\
  & +(-1)^{k}\phi(X_1,\cdots,X_k)X_{k+1}.$$
  \end{array}
$$ 
where $XY$ denotes the product on the algebra $A$. If this algebra is associative, then $\delta_H^{k+1} \circ \delta_H^k = 0.$

\medskip

Let us assume now that $(A, \ast)$ is a  \wa $\K$-algebra. To simplify notations, we write $XY$ in place of $X \ast Y$. For any $k \in \N^*$ let us consider the vector spaces
$$\mathcal{C}^k_{WA}(A)=\{\phi_k : A^{\otimes^k} \ra A\}$$
whose elements are the $k$-linear maps on $A$ with values in $A$. For $k=0$, we  put $\mathcal{C}^0_{WA}(A)=A.$ 
We shall define the first  coboundary operators $\delta_{WA}^k: \mathcal{C}^k_{WA}(A) \ra \mathcal{C}^{k+1}_{WA}(A).$
We put
$$\delta^0_{WA}(X)=L_X-R_X$$
where $L_X(Y)=XY$ and $R_X(Y)=YX$ for any  $X,Y \in A$ 
$$\delta^1_{WA}\phi_1(X,Y)=\phi_1(X)Y+X\phi_1(Y)-\phi_1(XY).$$
We verify that
$$\delta^1_{WA} \circ \delta^0_{WA} =0.$$
For the degree $2$ we put
\begin{equation}
\label{ delta2}
\delta_{WA}^2\phi_2(X_1,X_2,Z)=\delta_H^2\phi_2(X_1,X_2,Z)+\delta_H^2\phi_2(X_2,Z,X_1)-\delta_H^2\phi_2(X_2,X_1,Z)
\end{equation}
\medskip
or
$$\delta_{WA}^2\phi_2 = \delta_{H}^2\phi_2 \circ \Phi^A_{v_{WA}}$$
with $v_{WA}=Id+c_1-\tau_{12}.$ We verify easily the relation
$$\delta^2_{WA} \circ \delta^1_{WA} =0.$$

\medskip

The determination of the coboundary operators of greater degree is much more delicate. Focus on degree $3$. For this let us put
$$\begin{array}{l}
\delta^3_{WA}\varphi_3 (X_1,X_2,X_3,X_4)=\sum a_{i,jkl}X_i\varphi_3(X_j,X_k,X_l)+\sum b_{ijk,l}\varphi_3(X_i,X_j,X_k)X_l \\
+\sum c_{ij,k,l}\varphi_3 (X_iX_j,X_k,X_l)+\sum d_{i,jk,l}\varphi_3(X_i,X_jX_k,X_l)+\sum e_{i,j,kl}\varphi_3(X_i,X_j,X_kX_l).
\end{array}$$
the sums being on all permutations of $(1,2,3,4)$. This expression has to satisfy $\delta_{WA}^3 \circ \delta_{WA} ^2 = 0.$ But the reduction of this relation can be do using only the weakly associativity identity. But developing this identity, it appears three types of terms
\begin{enumerate}
\item terms of degree $3$: $X_i(X_j\varphi_2(X_k,X_l), \  (X_iX_j)\varphi_2(X_k,X_l),  \ X_i(\varphi_2(X_j,X_k)X_k)$, \\ $(X_i\varphi_2(X_j,X_k))X_l$ and the symmetric terms,
\item terms of degree $2$ $X_i\varphi_2(X_jX_k,X_l),X_i\varphi_2(X_j,X_kX_l)$ and the symmetric terms,
\item terms of degree $1$ $\varphi_2((X_iX_j)X_k,X_l),\varphi_2(X_i(X_jX_k),X_l),\varphi_2(X_i,(X_jX_k)X_l), \\ \varphi_2(X_i,X_j(X_kX_l)$
\item terms of degree $0$  $\varphi_2(X_iX_j,X_kX_l)$.
\end{enumerate}
Terms of degree $3$ will be reduced using the weakly associativity identity. Similarly, we must reduce the terms of degree 1 between them by using the same identity. However, no identity can reduce the terms of degree $2$ and degree $0$. They must therefore be eliminated algebraically. Let us consider the following vectors of $\KS$:
$$u_1=Id+c+c^2, \ u_2=\tau_{12}+c^2, \ u_3= \tau_{13}-c-c^2, \ u_4=\tau_{23}+c.$$
Thus, if $M=\{Id,(2134),(3124),(4123)\}\subset  \K[\Sigma_4]$ and if  $\Theta (X_1,X_2,X_3,X_4)$ refers $\delta_{WA}^3 \circ \delta_{WA} ^2 \varphi_2(X_1,X_2,X_3,X_4)$,  we have
$$
\begin{array}{l}
\displaystyle
\Theta(X_1,X_2,X_3,X_4)= \sum_{\sigma \in M }X_{\sigma(1)}\delta_H^2\varphi_2 \circ (a_{\sigma,1}u_1+a_{\sigma,2}u_2+a_{\sigma,3}u_3+a_{\sigma,4}u_4)(X_{\sigma(2)},X_{\sigma(3)},X_{\sigma(4)})\\
\displaystyle
+\sum_{\sigma \in M }(\delta_H^2\varphi_2 \circ (b_{\sigma,1}u_1+b_{\sigma,2}u_2+b_{\sigma,3}u_3+b_{\sigma,4}u_4)(X_{\sigma(2)},X_{\sigma(3)},X_{\sigma(4)})X_{\sigma(1)}\\
\displaystyle
+\sum_{\sigma \in\K[\Sigma_4]}\delta_H^2\varphi_2 \circ (c_{\sigma,1}u_1+c_{\sigma,2}u_2+c_{\sigma,3}u_3+c_{\sigma,4}u_4)(X_{\sigma(1)}X_{\sigma(2)},X_{\sigma(3)},X_{\sigma(4)}).\\
\end{array}
$$
This linear system,  $\delta_{WA}^3 \circ \delta_{WA} ^2 \varphi_2=0$  is a system wth $120$ variables and $360$ linear equations. Let us note also that if $(S)$ is a solution of this system, then $(S) \circ \Phi_v$ is also a solution for any $v \in \Sigma_4$. 
The vectors $u_1,u_2,u_3,u_4$ that appear in this linear system verify
 $u_1,u_3 \in \mathcal{O}(v_{WA})$ and $u_2,u_4 \in \K[\mathcal{O}(v_{WA})].$ Moreover we have $u_1,u_2,u_4 \in \K[\mathcal{O}(u_3)]$. This leads us to the identity
$$
\begin{array}{l}
\displaystyle
\delta_{WA}^3  \varphi_3(X_1,X_2,X_3,X_4)= \sum_{\sigma \in M }X_{\sigma(1)}\varphi_3 \circ (a_{\sigma,3}u_3)(X_{\sigma(2)},X_{\sigma(3)},X_{\sigma(4)})\\
\displaystyle
+\sum_{\sigma \in M }(\varphi_3 \circ b_{\sigma,3}u_3)(X_{\sigma(2)},X_{\sigma(3)},X_{\sigma(4)})X_{\sigma(1)}
+\sum_{\sigma \in\K[\Sigma_4]}\varphi_3 \circ (c_{\sigma,3}u_3)(X_{\sigma(1)}X_{\sigma(2)},X_{\sigma(3)},X_{\sigma(4)}).\\
\end{array}
$$

\subsection{The quadratic operad $\mathcal{W}ass$}

 We recall briefly some facts about operads. An operad $\mathcal{P}$ is a sequence of vector spaces $\mathcal{P}(n)$, for $n \geq  1$, such that $\mathcal{P}(n)$ is a module over the symmetric group $\Sigma_n$, together with composition maps
$$
\circ_i: \mathcal P(n) \times \mathcal P(m) \to \mathcal P(n+m-1)
$$
satisfying associativity-like conditions:
$$
(f \circ_i g) \circ_j h = f \circ_j (g \circ_{i-j+1} h) 
$$
 An algebra over an operad $\mathcal{P}$ is a vector space $A$ together with maps $\mathcal{P}(n) \otimes A^{\otimes^n} \ra A$ satisfying some relations of associativity, unitarity and equivariance.

The operad for weakly-associative 
algebras will be denoted by $\mathcal{W}ass$.
It is a binary quadratic operad, that is, operad of the form $\p =\Gamma(E)/(R)$, where $\Gamma(E)$ denotes the free operad generated by a
$\Sigma_2$-module $E$ placed in arity 2 and $(R)$ is the operadic ideal generated by a $\Sigma_3$-invariant subspace $R$
of $\Gamma(E)(3).$ To define $\mathcal{W}ass$, we consider the operadic ideal $R \subset \Gamma(E)(3)$ generated by the vector 
\begin{equation}
\label{R3}
(x_1x_2)x_3-x_1(x_2x_3)+(x_2x_3)x_1-x_2(x_3x_1)-(x_2x_1)x_3+x_2(x_1x_3).
\end{equation}
 From the study, in the previous sections of the orbit of the vector $Id+c_1-\tau_{12} \in \KS$, we deduce that the $\K$-vector space $R$ is of dimension $4$, generated by the vectors of $\Gamma(E)(3)$:
$$
 \left\{
 \begin{array}{ll}
(x_1x_2)x_3-x_1(x_2x_3)+(x_2x_3)x_1-x_2(x_3x_1)-(x_2x_1)x_3+x_2(x_1x_3)\\
(x_2x_1)x_3-x_2(x_1x_3)+(x_1x_3)x_2-x_1(x_3x_2)-(x_1x_2)x_3+x_1(x_3x_3)\\
(x_3x_2)x_1-x_3(x_2x_1)+(x_2x_1)x_3-x_2(x_1x_3)-(x_2x_3)x_1+x_2(x_3x_1)\\
(x_1x_3)x_2-x_1(x_3x_2)+(x_3x_2)x_1-x_3(x_2x_1)-(x_3x_1)x_2+x_3(x_1x_2)\\
\end{array}
  \right.
  $$
 Then the first terms of $\mathcal{W}ass=\oplus_{n \geq 1} \mathcal{W}ass(n)$ are
 $$\mathcal{W}ass(1)=\K,\mathcal{W}ass(2)=\Gamma(E)(2)=\K\{x_1x_2,x_2x_1\},\mathcal{W}ass(3)=\Gamma(E)(3)/R$$
 and $\dim \mathcal{W}ass(3)=8.$

Recall the definition of the quadratic dual operad $\p^! $ of an operad $\p$. If $\p = \Gamma(E)/(R)$, then there is a scalar product on $\Gamma(E)(3)$ defined by
\begin{eqnarray}
\label{pairing}
\left\{
\begin{array}{l}
<(x_i \cdot x_j)\cdot x_k,(x_{i'} \cdot x_{j'})\cdot x_{k'}>=0, \ {\rm if} \ \{i,j,k\}\neq \{i',j',k'\}, \\
<(x_i \cdot x_j)\cdot x_k,(x_i \cdot x_j)\cdot x_k>=(-1)^{\varepsilon(\sigma)},  \\
\qquad \qquad \qquad   {\rm with} \ \sigma =
\left(
\begin{array}{lll}
i&j&k \\
i'&j'&k'
\end{array}
\right)
\ {\rm if} \ \{i,j,k\}= \{ i',j',k'\}, \\
<x_i \cdot (x_j\cdot x_k),x_{i'} \cdot (x_{j'}\cdot x_{k'})>=0, \ {\rm if} \ \{ i,j,k\} \neq \{ i',j',k'\}, \\
<x_i \cdot (x_j\cdot x_k),x_i \cdot (x_j\cdot x_k)>=-(-1)^{\varepsilon(\sigma)}, \\
\qquad \qquad \qquad  {\rm with} \ \sigma =
\left(
\begin{array}{lll}
i&j&k \\
i'&j'&k'
\end{array}
\right)
\ {\rm if} \ (i,j,k)\neq (i',j',k'), \\
<(x_i \cdot x_j)\cdot x_k,x_{i'} \cdot (x_{j'}\cdot x_{k'})>=0. 
\end{array}
\right.
\end{eqnarray}
Let $R^\perp$ be the annihilator of $R$ with respect to this scalar product. If $R$ is the operatic ideal generated by (\ref{R3}), then $R^\perp$ is generated by the relations
$$
\left\{
\begin{array}{l}
(x_1x_2)x_3-x_1(x_2x_3), \\$$
(x_1x_2)x_3+(x_3x_2)x_1-(x_1x_3)x_2-(x_2x_3)x_1
\end{array}
\right.
$$
In fact, it is clear that $(x_1x_2)x_3-x_1(x_2x_3) \in R^\perp.$ Thus, it is sufficient to consider only vectors of type $(xy)z$. Let be
$$u=a(x_1x_2)x_3+b(x_2x_1)x_3+c(x_3x_2)x_1+d(x_1x_3)x_2+e(x_2x_3)x_1+f(x_3x_1)x_2.$$
If we write that $u \in R^\perp$, we obtain
 $$
 \left\{
 \begin{array}{ll}
 a+b+e=0\\
 -b-d-a=0\\
 -c-b-e=0\\
 -d-c-f=0.
\end{array}
  \right.
  $$
that is 
$$e=d=-a-b,c=a,f=b$$
and
 $$u=a((x_1x_2)x_3+(x_3x_2)x_1-(x_1x_3)x_2-(x_2x_3)x_1)+b((x_2x_1)x_3-(x_1x_3)x_2-(x_2x_3)x_1+(x_3x_1)x_2).$$
 This gives the generator of the operatic ideal $R^\perp$. Since the dual operad $\mathcal{W}ass^!$ is defined by the operatic ideal $R^\perp$, then
 $$\dim \mathcal{W}ass^! (3)   =4.$$

\begin{propo}
A $\mathcal{W}ass^!$-algebra is an associative algebra  satisfying the relation
$$abc+cba-acb-bca.$$
\end{propo}

Let us consider now $\mathcal{W}ass^!(4)$. It can be considered as a linear subspace of $\K[\Sigma_4]$ defined by relations coming from $R(3)$ and concerning the vectors of type $$(\bullet \bullet)\bullet\bullet,\bullet(\bullet\bullet)\bullet,\bullet\bullet(\bullet\bullet),(\bullet\bullet\bullet)\bullet,\bullet(\bullet\bullet\bullet).$$
Let us denote by $\mathcal{R}_3(a,b,c)$ the relation $abc+cba-acb-bca.$ We have
$$
\left\{
\begin{array}{l}
\mathcal{R}_3(a,b,c)+\mathcal{R}_3(a,c,b)=0,\\
\mathcal{R}_3(a,b,c)+\mathcal{R}_3(c,a,b)-\mathcal{R}_3(b,a,c)=0.\\
\mathcal{R}_3(a,bc,d)+\mathcal{R}_3(a,db,c)-\mathcal{R}_3(a,cd,b)=0, \\
\mathcal{R}_3(a,b,c)d-\mathcal{R}_3(d,bc,a)+\mathcal{R}_3(d,cb,a)-d\mathcal{R}_3(a,b,c)=0.
\end{array}
\right.
$$
We deduce that the relations $\mathcal{R}_3(ab,c,d)=0$ are consequences of the relations $\mathcal{R}_3(d,ab,c)=\mathcal{R}_3(c,ab,d)=0$
We have also seen that the vector
$$v=x_1x_2x_3+x_3x_2x_1-x_1x_3x_2-x_2x_3x_1$$
satisfies $\tau_{23}(v)=-v.$ This implies that the relations concerning $\bullet (\bullet \bullet)\bullet $ and $\bullet\bullet(\bullet\bullet)$ are identical.  The second relation shows that any relation of type $\mathcal{R}_3(\bullet\bullet,\bullet,\bullet)$ is a linear combination of relations of type $\mathcal{R}_3(\bullet,\bullet\bullet,\bullet)$. The third relation shows that any relation of type $\mathcal{R}_3(\bullet,\bullet\bullet,\bullet)$ is a linear combination of the $16$ relations
$$
\left\{
\begin{array}{llll}
\mathcal{R}_3(1,23,4) & \mathcal{R}_3(2,13,4) & \mathcal{R}_3(3,21,4) & \mathcal{R}_3(4,23,1)\\
\mathcal{R}_3(1,32,4) & \mathcal{R}_3(1,43,2) & \mathcal{R}_3(2,31,4) & \mathcal{R}_3(2,43,1)\\
\mathcal{R}_3(3,12,4) & \mathcal{R}_3(3,24,1) & \mathcal{R}_3(4,13,2) & \mathcal{R}_3(4,21,3)\\
\mathcal{R}_3(1,34,2) & \mathcal{R}_3(2,34,1) & \mathcal{R}_3(3,42,1) & \mathcal{R}_3(4,31,2)
\end{array}
\right.
$$

The fourth relation shows that $\mathcal{R}_3(\bullet,\bullet,\bullet)\bullet$ \ is a linear combination of relations of type
$\mathcal{R}_3(\bullet,\bullet\bullet,\bullet)$ and $\bullet\mathcal{R}_3(\bullet,\bullet,\bullet).$ We deduce that we have only to consider relations of type 
$\mathcal{R}_3(\bullet,\bullet\bullet,\bullet)$ and $\bullet\mathcal{R}_3(\bullet,\bullet,\bullet).$

\begin{lemma}
We have the relation
$$
\begin{array}{l}
\mathcal{R}_3(1,43,2)-\mathcal{R}_3(1,34,2) +\mathcal{R}_3(2,34,1)-\mathcal{R}_3(2,43,1)+\mathcal{R}_3(3,21,4)-\mathcal{R}_3(3,12,4)\\
+\mathcal{R}_3(4,12,3)-\mathcal{R}_3(4,21,3)=0
\end{array}
$$
\end{lemma}
We deduce
\begin{propo}
$$\dim \mathcal{W}ass^! (4)   =6.$$
\end{propo}
\pf In fact, if we consider the matrix corresponding to the relations $\mathcal{R}_3(\bullet,\bullet\bullet,\bullet)=0$ and $\bullet\mathcal{R}_3(\bullet,\bullet,\bullet)=0$, it is a square matrix and, from the previous remark of rank equal to $18$. Its kernel is of dimension $6$. The relations of definition of $\mathcal{W}ass^! (4)$ are 
\begin{equation}
\label{W4dual}
\left\{
\begin{array}{l}
x_1x_2x_3x_4+x_4x_2x_3x_1-x_1x_4x_2x_3-x_2x_3x_4x_1=0 \\
x_1x_2x_3x_4+x_1x_4x_3x_2-x_1x_2x_4x_3-x_1x_3x_4x_2=0
\end{array}
\right.
\end{equation}

\noindent{\bf Remark.} Recall that the generating function of a quadratic operad $\mathcal{P}=\oplus \mathcal{P}(n)$ is the formal series
$$f_\mathcal{P}(t)=\sum_{i=1}^{+ \infty}(-1)^n \dim \mathcal{P}(n)\frac{x^n}{n!}.$$
Then
$$f_{ \mathcal{W}ass^!}=-x+2\frac{x^2}{2}-4\frac{x^3}{6}+6\frac{x^4}{24}- \cdots=-x+x^2-\frac{2x^3}{3}+\frac{x^4}{4}- \cdots$$
and
$$f_{ \mathcal{W}ass}=-x+2\frac{x^2}{2}-8\frac{x^3}{6}+ \cdots =-x +x ^2-\frac{4x^3}{3}+ \cdots .$$

\medskip

\noindent{\bf The operadic cohomology.} Let $\mathcal{P} =\oplus \mathcal{P}(n)$ be a quadratic operad and assume
$A $ is a $\mathcal{P}$-algebra.  An operadic cochain complex is defined by:
$$\mathcal{C}^k(A,A)=Hom (\mathcal{P}^{!})^{\vee}(k)\otimes_{\Sigma_n} A^{\otimes ^k},A)$$
where $(\mathcal{P}^{!})^{\vee}=Hom(\mathcal{P}^{!}),\K) \otimes Sgn_n.$
In particular $$\mathcal{C}^2(A,A)=Hom (A^{\otimes ^2},A).$$
\begin{definition}A $3$-cochain for the operadic complex is a trilinear map $\Psi$ on $A$ which satisfies 
$$\Psi(x_1,x_2,x_3)+\Psi(x_2,x_1,x_3)-\Psi(x_1,x_3,x_2)-\Psi(x_2,x_3,x_1)=0.$$
\end{definition}
 In fact $\mathcal{W}ass^{!}(3)$ is generated by the relation
$$x_1x_2x_3+x_3x_2x_1-x_1x_3x_2-x_2x_3x_1=0$$
which corresponds to the vector $Id+\tau_{13}-\tau_{23}-c_1$. 
Let $w \in \KS$ and let us consider the equation $(Id-\tau_{12}+c)\circ w= 0$. Then $w$ is a linear combination of the vectors
$$w_1=Id+\tau_{12}-\tau_{23}-c, \ w_2=\tau_{12}-\tau_{13}-c+c^2.$$
We verify that 
$$\delta^2_{WA}\varphi \circ \Phi^A_{w_1}=0$$
this justifies the previous definition.

\medskip

Concerning  degree $4$, we have seen that the relation of definition of $\mathcal{W}ass^{!}(4)$ given in (\ref{W4dual}) are
$$
\left\{
\begin{array}{l}
x_1x_2x_3x_4+x_4x_2x_3x_1-x_1x_4x_2x_3-x_2x_3x_4x_1=0 \\
x_1x_2x_3x_4+x_1x_4x_3x_2-x_1x_2x_4x_3-x_1x_3x_4x_2=0
\end{array}
\right.
$$
and correspond to the vectors
$$
\left\{
\begin{array}{l}
Id+(14)-(243)-(1234)=0 \\
Id+(24)-(34)-(234)=0
\end{array}
\right.
$$
where $(i_1i_2,\cdots,i_n)(j)=j$ if $j \neq i_k, k+1,\cdots,n$ and $(i_1i_2,\cdots,i_n)(j)=i_{k+1} \ (\rm{mod} \ n )$ if $j=i_k$.
The action of $\K[\Sigma_4]$ on $\mathcal{W}ass^{!}(4)$ corresponds to the vectors
$$
\left\{
\begin{array}{l}
Id+(14)-(234)-(1432)=0 \\
Id+(24)-(34)-(243)=0
\end{array}
\right.
$$
\begin{proposition}
A $4$-cochain for the operadic complex is a $4$-linear map $\Theta$ on $A$ which satisfies 
$$
\left\{
\begin{array}{l}
\Theta(x_1,x_2,x_3,x_4)+\Theta(x_4,x_2,x_3,x_1)-\Theta(x_1,x_3,x_4,x_2)-\Theta(x_4,x_1,x_2,x_3)=0 \\
\Theta(x_1,x_2,x_3,x_4)+\Theta(x_1,x_4,x_3,x_2)-\Theta(x_1,x_2,x_4,x_3)-\Theta(x_1,x_4,x_3,x_2)=0. \\
\end{array}
\right.
$$
\end{proposition}
We deduce that, if $\Psi \in  \mathcal{C}^3(A,A)$, then

$\delta^3_{WA}\Psi(x_1,x_2,x_3,x_4)+\delta^3_{WA}\Psi(x_4,x_2,x_3,x_1)-\delta^3_{WA}\Psi(x_1,x_3,x_4,x_2)-\delta^3_{WA}\Psi(x_4,x_1,x_2,x_3)=0$

$\delta^3_{WA}\Psi(x_1,x_2,x_3,x_4)+\delta^3_{WA}\Psi(x_1,x_4,x_3,x_2)-\delta^3_{WA}\Psi(x_1,x_2,x_4,x_3)-\delta^3_{WA}\Psi(x_1,x_4,x_3,x_2)=0.$

\medskip

\noindent{\bf Consequence.} Let be $\varphi \in \mathcal{C}^2(A,A)$. Then $\delta^2_{WA}\varphi  \in  \mathcal{C}^3(A,A)$ and if $v_4=Id+(14)-(234)-(1432)=0$, 
$v'_4=Id+(24)-(34)-(243)=0$, then 
$$
\left\{
\begin{array}{l}
\delta^3_{WA}(\delta^2_{WA}\varphi )\circ \Phi^A_{v_4}=0,\\
\delta^3_{WA}( \delta^2_{WA}\varphi )\circ \Phi^A_{v'_4}=0. \\
\end{array}
\right.
$$

\section{quantization of nonassociative Poisson algebras}
\subsection{Recall: associative case}
Let $(A,\mu)$ be an associative algebra.  A formal deformation of $(A,\mu)$ is given by a family 
$$\{\varphi_i : A \otimes A \rightarrow A , \  i \in \N\}$$
with $\mu_0=\mu$ and
satisfying
\begin{equation}
\label{def}
\sum_{i+j=k, i,j \geq 0}\varphi_i(\varphi_j(X,Y),Z)=\sum_{i+j=k, i,j \geq 0}\varphi_i(X,\varphi_j(Y,Z))
\end{equation}
for any $X,Y,Z \in A$ and $k \geq 1.$ If we denote by $\K[[t]]$ the algebra of formal series with one indeterminate $t$, this definition is equivalent  to consider on the space $A [[t]] $ of formal series with coefficients in $A$ a structure $\mu_t$ of $\K[[t]]$-associative algebra such that the canonical map $A[[t]]/\K[[t]] \rightarrow A$ is an isomorphism of algebras. It is practical to write
$$\mu_t=\mu+t \varphi_1+t^2\varphi_2+ \cdots$$
Equation (\ref{def}) implies at the order $k=0$ that $\mu$ is associative, at the order $k=1$ that $\varphi_1$ is a $2$-cocycle for the Hochschild cohomology of $(A,\mu)$, that is
$$\delta_H^2\varphi_1(X,Y,Z)=X\varphi_1(X,Y)-\varphi(XY,Z)+\varphi_1(X,YZ)-\varphi_1(X,Y)Z$$
where, to simplify the notations, $XY$ means $\mu(X,Y)$.
\begin{lemma}
Let $\psi$ be a skew-symmetric bilinear map on $A$. Then the Leibniz identity
$$\psi(XY,Z)-X\psi(Y,Z)-\psi(X,Z)Y=0$$
is equivalent to
$$\delta_H^2\psi =0.$$
\end{lemma}
\pf In a first step, let us show that the Leibniz identity is equivalent to
$$\psi(XX,Y)=2X\psi(X,Y)$$
for any $X,Y \in A.$ In fact, assume that $\psi(XY,Z)-X\psi(Y,Z)-\psi(X,Z)Y=0.$
 Then, for $Y=X$ we obtain
 $$\psi(XX,Z)=X\psi(X,Z)+\psi(X,Z)X
$$
and since $\mu(X,Y)=XY$ is commutative,
 $$\psi(XX,Z)=2X\psi(X,Z)$$
 Conversely, if for any $X,Z \in A$ we have $\psi(XX,Z)=2X\psi(X,Z)$ , then by
  linearization, we obtain
 $$2\psi(XY,Z)=2Y\psi(X,Z)+2X\psi(Y,Z)$$
 that is the Leibniz identity.
 Let us come back to the proof of the lemma. Assume that $\psi$ satisfies the Leibniz identity on $A$:
$$\psi(XY,Z)=X\psi(Y,Z)+\psi(X,Z)Y.$$
This gives
$$
\begin{array}{ll}
    \delta_H^2 \psi (X,Y,Z)  & = X\psi(Y,Z)-\psi(XY,Z)+\psi(X,YZ)-\psi(X,Y)Z   \\
      &   = X\psi(Y,Z)-X\psi(Y,Z)-\psi(X,Z)Y+Y\psi(X,Z)+\psi(X,Y)Z-\psi(X,Y)Z\\
      &=0.
\end{array}
$$
Conversely, if $\psi$ is a skew-symmetric $2$-cocycle of $(A,\mu)$, then
$$X\psi(Y,Z)-\psi(XY,Z)+\psi(X,YZ)-\psi(X,Y)Z=0.$$
In particular, for $Y=X$, 
$$X\psi(X,Z)-\psi(XX,Z)+\psi(X,XZ)-\psi(X,X)Z=X\psi(X,Z)-\psi(XX,Z)+\psi(X,XZ)=0$$
and if $Z=X$:
$$X\psi(Y,X)-\psi(XY,X)+\psi(X,YX)-\psi(X,Y)X=2X\psi(Y,X)+2\psi(X,XY)=0.$$
We deduce
$$\psi(XX,Z)=X\psi(X,Z)+\psi(X,XZ)=2\psi(X,XZ)$$
and $\psi$ satisfies the Leibniz identity.
\begin{lemma}
Let $\varphi_1$ be the term of degree $1$ of a formal deformation of the associative commutative product $\mu$. Then $\varphi_1$ is Lie admissible.
\end{lemma}
\pf Equation (\ref{def}) at the order $2$ gives
$$\varphi_1(\varphi_1(X,Y),Z)-\varphi_1(X,\varphi_1(Y,Z)=\delta_H^2\varphi_2(X,Y,Z).$$
But
$$\delta_H^2\varphi_2 \circ \Phi^A_{Id-\tau_{12}-\tau_{13}-\tau_{23}+c+c^2}=0$$
this implying
$$\mathcal{A}(\varphi_1)\circ \Phi^A_{Id-\tau_{12}-\tau_{13}-\tau_{23}+c+c^2}=0$$
this means that $\varphi_1$ is Lie admissible. Since $\delta \varphi_1=0$ implies, because $\mu$ is commutative, that $\delta^2_H\psi_{\varphi_1}=0$, 
we have:
\begin{proposition} The skew-symmetric part $\psi_{\varphi_1}$ of $\varphi_1$ is a Lie bracket on $A$ and $(A,\mu,\psi_{\varphi_1})$ is a Poisson algebra. The formal deformation $(A[[t]],\mu_t)$ is a quantization of the Poisson algebra $(A,\mu,\psi_{\varphi_1})$. 
\end{proposition}
Now, the problem is to know if any Poisson algebra $(A,\mu,\psi)$ come from a formal deformation of $(A,\mu)$
 
 \medskip
 
 \noindent{\bf Case of finite dimension.}  We assume that $\K$ is algebraically close. Let $Ass_n$ be the algebraic variety whose elements are the $n$-dimensional associative multiplication on $\K^n$. This variety is fibered by the orbits of the natural action of the algebraic group $GL(n,\K)$, the orbit of $\mu$ being the set of isomorphic multiplication to $\mu$. A multiplication $\mu \in Ass_n$ is called rigid if its orbit is Zariski open in $Ass_n$. For exemple, if the second group of the Hochschild cohomology of $\mu$ is equal to $\{0\}$, then $\mu$ is rigid. But there exists also rigid multiplication with a non trivial cohomology group. In this case, if $\mathcal{Ass}_n$ denotes the associated affine schema, this schema is not reduced. In particular, for a rigid multiplication, there exists a $2$-cocycle $\varphi$ which is not integrable, that is it is not the first term of a formal deformation of $\mu$. If its skew-symmetric part is a Lie bracket, then we obtain in this way a Poisson structure on $(\K^n,\mu)$ not associated with a deformation. For example, let us consider the associative commutative multiplication with an idempotent $e$ such that
 $$e.e_i=e_i, i=1,\cdots,p, \ \ ef_j=0, \ j=1,\cdots,q$$
 with $p+q+1=n$ and $\{e,e_i,f_j\}$ a basis of $\K^n$. Since $e^2=e$, the associativity conditions imply that we have the direct decomposition
 $\K^n= A \oplus B$
 with $A=\K\{e,e_i\}$, $B=\K\{f_j\}$ and $AB=\{0\}$.  The subspace $A$ is an unitary commutative subalgebra of $\K^n$ and $B$ a commutative subalgebra. If $\{,\}$ is a Poisson bracket on $(\K^n,\mu)$, then
 $$\{e^2,e_i\}=\{e,e_i\}=2e\{e,e_i\}=0$$
 because $1/2$ is not an eigenvalue of the linear map $L_e: x \rightarrow ex.$ Likewise $\{e,f_j\}=0.$  The computations of $\{e_i,ef_j\}$ and $\{f_j,ee_i\}$ imply that $\{e_i,f_j\}=0.$ We deduce that for any $2$-cocycle $\varphi$, its skew-symmetric part satisfies
 $$\psi_\varphi(e,A\oplus B)=0, \psi_\varphi(A,B)=0.$$
 If $A$ is a rigid unitary commutative associative algebra with $H^2_H(A,A)\neq 0$, then a non symmetric non integrable cocycle   determines a Poisson algebra which is not associated with a quantization.
 
 \medskip
 
 \noindent{\bf Skew-symmetric multi-derivations.}  A skew-symmetric map on $(A,\mu)$ is a biderivation if it is a derivation on each argument, that is if it satisfies the Leibniz identity. A skew-symmetric $n$-linear map is a $n$-derivation (or multi-derivation in general)
 if it is a derivation on each argument. Lemma 8 can be interpreted saying that the skew-symmetric map $\psi$ is a biderivation if and only if it is a skew-symmetric $2$-cocycle. This result can be generalized. For example, if $\psi$ is a skew-symmetric $3$-linear map. Then if it is a $3$-cocycle, we have
 $$X\psi(Y,Z,T)-\psi(XY,Z,T)+\psi(X,YZ,T)-\psi(X,Y,ZT)-\psi(X,Y,Z)T=0.$$
 In particular, we have
 $$
 \left\{
 \begin{array}{l}
 \psi(X^2,Y,Z)=2X\psi(X,Y,Z),\\
  \psi(X,Y,XZ)=X\psi(X,Y,Z),\\
    \psi(X,Y,XZ)=\psi(X,XY,Z).
    \end{array}
    \right.
    $$
If we linearize these identities, we obtain
 $$
 \left\{
 \begin{array}{l}
 \psi(XT,Y,Z)=T\psi(X,Y,Z)+X\psi(T,Y,Z), \\
 \psi(X,Y,ZT)=X\psi(T,Y,Z)+T\psi(X,Y,Z),\\
 \psi(X,Y,TZ)+\psi(T,Y,XZ)=\psi(T,XY,Z)+\psi(X,TY,Z).
  \end{array}
    \right.
    $$
This shows that $\psi$ is a $3$-derivation.  Conversely, if $\psi$ is a $3$-linear skew-symmetric $3$-derivation, that is
$$ \psi(XT,Y,Z)=T\psi(X,Y,Z)+X\psi(T,Y,Z)$$
for any $X,Y,Z,T \in A.$ 
%Then we have
%$$\left\{
% \begin{array}{l}\psi(X^2,Y,Z)=2X\psi(X,Y,Z),\\
% \psi(Y,X,XZ)=\psi(XZ,Y,X)=X\psi(X,Y,Z),\\
%  \psi(X,Y,XZ)=-\psi(XZ,Y,X)=-X\psi(Z,Y,X)=X\psi(Y,Z,X)=\psi(XY,Z,X).  \end{array}
 %   \right.
  %  $$
 A direct computation shows that this identity implies that the Leibniz identity is satisfied. Then for any skew-symmetric $3$-ary applications $\psi$, the Leibniz identity on $(A,\mu)$ is equivalent to $\delta_H^3(\psi)=0.$ And this is true for the greater order.

\subsection{Deformation of weakly associative algebras}

\begin{definition}( Deformation of weakly associative algebras.) Let $(A,\mu)$ be a weakly associative $\K$-algebra. 
A deformation $\mu_t$ of $\mu$ is a bilinear map
%\footnote{The map $C$ can be straightforwardly linearly extended to a map $C:V\eps\otimes V\eps\to V\eps$. } 
$\mu_t: A\otimes A\ra A[[t]]$ of the form:
$$
\mu_t=\mu+t\, \varphi_1+t^2\, \varphi_2+\ldots \label{eqexpansion}
 $$
 which is weakly associative 
 where each $\varphi_i$ is a bilinear map $\varphi_i:A\otimes A\to A$. 
 
 The map $\mu_t$ can be straightforwardly extended linearly to a map $\mu_t:A[[t]]\otimes A[[t]]\to A[[t]]$, thus endowing the vector space $A[[t]]$ with a structure of weakly associative $\K[[t]]$-algebra. 
\end{definition}

\begin{definition} (Equivalence of deformations of weakly associative algebras). Let $(A,\mu)$ be a weakly associative $\K$-algebra. 
Two deformations $\mu_t,\mu'_t$ will be said equivalent if there exists an isomorphism of $\K[[t]]$-algebras:
$$
\phi:(A[[t]],\mu_t)\overset{\sim}{\to}(A[[t]],\mu'_t)
$$
such that $\phi(X)=X+\mathcal O(t)$ for all $X\in A$. 
\end{definition}
%\Rem{This condition implies isomoprhism of algebras but is stronger, as we ask that the zeroth order coincides with the identity, \ie $\varphi$ is invisible at the classical level.}
%\bdefi{Equivalence of deformations}{\label{defieq}Let $( \alg,\mu)$ be an associative $\mathbb R$-algebra and $t$ a formal parameter. Two formal deformations $\pset{B'_n}$ and $\pset{B_n}$ will be said equivalent if the respective associative $\mathbb Rt$-algebras $( \algt,*')$ and $( \algt,*)$ are isomorphic. }
% $\pset{B'_n}$ and $\pset{B_n}$ namely when the respective associative $\mathbb Rt$-algebras $( \algt,*')$ and $( \algt,*)$ are isomorphic. 
%Let us denote $A:( \algt,*')\to( \algt,*)$ such an isomorphism. The latter 

The previous condition is an equivalence relation and thus allows to define equivalence classes of deformations (denoted $[\mu]$). Two equivalent deformations  will sometimes be said gauge equivalent and acting on a deformation with a morphism $f_t$ will be referred to as a gauge transformation. We will call trivial a deformation gauge-equivalent to the zero algebra $(A[[t]],0)$.  

Explicitly, the above map $\phi$ must satisfy the two following conditions:
\begin{enumerate}
\item $f_t$ is an automorphism of the $\K[[t]]$-vector space $A[[t]]$ whose zeroth order coincides with the identity. 
\item $f_t$ is a homomorphism of $\K[[t]]$-algebras. 
\end{enumerate}
The first condition ensures that $f_t$ can be expanded as:
$$
f_t:=Id+\sum_{n=1}^\infty t^n\, h_n\label{eqvarphi}
$$
where $Id$ is the identity isomorphism on $A$ and $h_n$ are linear endomorphisms of $A$. 

%\Rem{Can 1 be replaced by any automorphism of $(V,\mu)$?}

The second condition reads explicitly as $f_t(\mu_t(X,Y))=\mu'_t(f_t(X), f_t(Y))$. It is clear that 
 $\mu'_t$ is weakly associative if and only if $\mu_t$ is.
 
\subsection{Deformation and polarization}
 
Let $(A,\mu)$ be a  weakly associative $\K$-algebra and $\mu_t=\mu+\sum t^i\varphi_i$ a weakly associative formal deformation of $\mu$. We denote by
$$\{X,Y\}=\mu(X,Y)-\mu(Y,X), \ \ X \bullet Y =\mu(X,Y)+\mu(Y,X)$$
the polarized multiplications associated with $\mu$ and
$$\{X,Y\}_t=\mu_t(X,Y)-\mu_t(Y,X), \ \ X \bullet_t Y =\mu_t(X,Y)+\mu_t(Y,X)$$
the polarized multiplications associated with $\mu_t.$ We have seen that $\{.,.\}$ and $\{.,.\}_t$ are Lie brackets respectively on $A$ and $A[[t]]$, $.\bullet .$ and $. \bullet_t . $ are commutative multiplications on these algebras and $(A,\bullet,\{.,.\})$ and $(A[[t]],\bullet_t,\{.,.\}_t)$ are respectively nonassociative Poisson $\K$-algebra and nonassociative Poisson $\K[[t]]$-algebra. Moreover $\{.,.\}_t$ is a Lie deformation of $\{.,.\}$ and $\bullet_t$ a commutative deformation of $\bullet$. Let us put
$$\{.,.\}_t=\{.,.\}+t B_1 + t^2 B_2 + \cdots $$
with $B_i(X,Y)=\varphi_i(X,Y)-\varphi_i(Y,X)$ and
$$\bullet_t=\bullet +t\rho_1+t^2\rho_2+\cdots$$
with $\rho_i(X,Y)=\varphi_i(X,Y)+\varphi_i(Y,X)$ for any $X,Y \in A.$  Since $\{.,.\}_t$ is a Lie deformation of the Lie bracket $\{.,.\}$, then the bilinear map $B_1$ is a $2$-cocycle for the Chevalley-Eilenberg cohomology of the Lie algebra $(A,\{.,.\})$  whose Jacobi expression $J(B_1)$ is a $2$-coboundary for this cohomology. In this case we shall say that $B_1$ is a Lie bracket modulo a coboundary.  In a same way, $\bullet_t$ is a commutative deformation of the commutative multiplication $\bullet$ and each bilinear map $\rho_i$ is symmetric. Let us consider now the $t$-expansion of the Leibniz identity
$$\{X,Y\bullet_t Z\}_t-Y\bullet_t \{X,Z\}_t-\{X,Y\}_t \bullet_t Z=0.$$
We obtain
\begin{enumerate}[(a)]
  \item ordre $0$:
$$
\{X,Y\bullet Z\}-Y\bullet \{X,Z\}-\{X,Y\} \bullet Z=0
$$ and we find again that $(A,\bullet,\{.,.\})$ is a nonassociative Poisson algebra,
  \item ordre $1$: 
$$
\{X,\rho_1(Y,Z)\}-\rho_1(Y,\{X,Z\})-\rho_1(\{X,Y\},Z)=-B_1(X,Y\bullet Z)+Y\bullet B_1(X,Z)+B_1(X,Y)\bullet Z.
$$
\end{enumerate}

\subsection{quantization by deformation}
Let $(A,\mu)$ be a commutative weakly associative $\K$-algebra. As previously, we write $XY$ in place of $\mu(X,Y)$. Let  $\mu_t = \mu + t \varphi_1+ t^2 \varphi_2 + \cdots $ be a weakly associative deformation of $\mu$. 
Writing that  $\mu_t$ is a weakly associative multiplication, the development at the order $1$ gives  $\delta_{WA}^2 \varphi_1=0 $ and at  the order $2$ the condition
$$\varphi_1 \circ (\varphi_1  \otimes Id)-Id \otimes \varphi_1) \circ \Phi^A_{Id+c_1-\tau_{12}}=\mathcal{WA}(\varphi_1)=\delta_{WA}^2 \varphi_2.$$
But for any  linear map $\varphi$, we have
$$\delta_{WA}^2 \varphi \circ \Phi^1_{Id-\tau_{23}}=0.$$
This implies
$$\mathcal{WA}(\varphi_1)\circ \Phi^1_{Id-\tau_{23}}=\mathcal{A}(\varphi_1)\circ \Phi^1_{v_{LAd}}=\delta_{WA}^2 \varphi_2 \circ \Phi^1_{Id-\tau_{23}}=0$$
and $\varphi_1$ is a Lie admissible multiplication.
\begin{proposition}
Let $\mu_t = \mu + t \varphi_1+ t^2 \varphi_2 + \cdots $ be a weakly associative deformation of the commutative weakly associative multiplication $\mu$. Then the bilinear map $\varphi_1$ is a $2$-cocycle for the $WA$-cohomology of $\mu$, that is 
$$\delta_{WA}^2 \varphi _1 =0$$
and it is Lie admissible, that is its skew-symmetric part $\psi_{\varphi_1}$ is a Lie bracket.
\end{proposition}

\medskip

%%%%%%%%%%%%%
%%%%%%%%%%%%%%%
%%%%%%%%%%%%%%%

\medskip

\begin{lemma}
Let $\varphi_1$ be a $2$-cocycle for the $WA$-cohomology of the weakly associative commutative algebra $(A,\mu)$. Then its skew-symmetric part $\psi_{\varphi_1}$ satisfies the Leibniz identity with respect to the weakly associative product:
$$\psi_{\varphi_1}(XY,Z)-X\psi_{\varphi_1}(Y,Z)-\psi_{\varphi_1}(X,Z)Y=0$$
for any $X,Y,Z \in A$ where $XY=\mu(X,Y)$ is the  weakly associative multiplication.
\end{lemma}
\pf In fact, since $\mu$ is commutative, we have
$$\delta_{WA}^2 \varphi (X,Y,Z)=\psi_{\varphi}(X,YZ)-Y\psi_{\varphi}(X,Z)-\psi_{\varphi}(X,Y)Z.$$
If $\delta_{WA}^2 \varphi=0$, then $\psi_{\varphi}$ satisfies the Leibniz identity. 

We can summarize the previous result:
\begin{proposition}
Let $(A,\mu)$ be a commutative weakly associative $\K$-algebra. Any weakly associative deformation 
$$\mu_t = \mu + t \varphi_1+ t^2 \varphi_2 + \cdots $$
of $\mu$ as a weakly associative algebra determines  a nonassociative Poisson algebra $(A,\mu,\psi_{\varphi_1})$
The formal deformation $(A[[t]],\mu_t)$ is called a quantization of the nonassociative Poisson algebra $(A,\mu,\psi_{\varphi_1})$.
\end{proposition}

\medskip

\noindent{\bf Remark.} Let $\mu(X,Y)=XY$ be a commutative multiplication, without other hypothesis. For any skew- symmetric bilinear map $\psi$, we can consider the operators
$$\delta^2_H \psi (X,Y,Z)=X\psi(Y,Z)-\psi (XY,Z)+ \psi (X,YZ)- \psi (X,Y) Z$$
and
$$\delta_{WA}^2 \psi = \delta^2_H \psi \circ \Phi^A_{Id+c_1-\tau_{12}}.$$
If $\mathcal{L}(\psi)$ denotes the Leibniz identity with respect to $\mu$:
$$\mathcal{L}(\psi)(X,Y,Z)=\psi(XY,Z)-X\psi(Y,Z)-\psi (X,Z)Y$$
then we have seen that the commutativity of $\mu$ and the skew-symmetry property of $\psi$ imply
$$\delta_{WA}^2 \psi  (X,Y,Z)=-2\mathcal{L}(\psi) (Y,Z,X)$$
and
$$2 \delta_H^2 \psi(X,Y,Z)=-\mathcal{L}(\psi)(X,Y,Z)-\mathcal{L}(\psi)(Y,Z,X).$$ 
Since $\delta_H^2 (\psi)=0$ implies $\delta^2_{WA}\psi=0$,
we deduce 
\begin{equation}
\label{l}
\mathcal{L}(\psi)=0 \Leftrightarrow \delta_{WA}^2 \psi=0 \Leftrightarrow \delta_{H}^2 \psi=0.
\end{equation}
Let us remark also that if $\rho$ is a symmetric bilinear map, then $\delta_{WA}^2\rho=0$ and $\delta_H^2 \rho \circ \Phi^A_{Id+\tau_{13}}=0.$ In particular, if $(A,\mu)$ is a commutative weakly associative algebra, then $\delta^2_{WA} \varphi =0$ if and only if 
$\delta^2_{WA}\psi_{\varphi}=0.$

\medskip

\noindent {\bf Applications: quantization of associative commutative algebras}
Considering a commutative associative algebra $(A,\mu)$, we have recalled in the first part how any $2$-cocycle of the Hochschild cohomology of $A$ which is the first term of a formal deformation of $\mu$ permits to construct on $A$ a Poisson structure and this deformation is a quantization of this Poisson algebra.
The previous result extends this result considering deformation of $\mu$ in the class of weakly associative algebra.
\begin{proposition}
Let  $(A,\mu)$ is an associative commutative algebra. Then any weakly associative formal deformation $\mu_t$ of $\mu$  defines a classical Poisson algebra $(A,\mu,\psi_{\varphi_1})$ and the algebra $(A[[t]],\mu_t)$  is a quantization of the Poisson algebra $(A,\mu,\psi_{\varphi_1})$.  The Poisson structure  so constructed on the algebra $A$ only depends of the equivalence class of formal deformations of $\mu$.
\end{proposition}
%For example, let us consider the $2$-dimensional associative commutative algebra given by
%$$\left\{\begin{array}{l}\mu(e_1,e_1)=e_1, \\\mu(e_1,e_2)=\mu(e_2,e_1)=\mu(e_,e_2)=0.\end{array}\right.$$Let us consider the formal deformation$$
%\left\{
%\begin{array}{l}\mu(e_1,e_1)=0, \\\mu(e_1,e_2)=\epsilon e_1,\\\mu(e_2,e_1)=2\epsilon e_1,\\\mu(e_,e_2)=3\epsilon e_2\end{array}\right.$$
%This algebra is not associative but weakly associative. It defines quantization of the Poisson algebra whose Poisson bracket is given by $$\{e_1,e_2\}_1=-e_1.$$
A special class of weakly associative formal deformation is given by the linear deformation: $\mu_t=\mu + t \psi_1$ with $\psi_1$ is a skew symmetric bilinear map such that $\delta_{WA}^2 \psi_1=0$ and $\mathcal{WA}(\psi_1)=0.$  Since $\psi_1$ is skew symmetric, the condition $\mathcal{WA}(\psi_1)=0$ is equivalent to say that $\psi_1$ is a Lie bracket. Then $\psi_1$ defines a Poisson structure on the associative commutative algebra $(A,\mu)$ associated with the linear quantization $(A[[t]],\mu_t=\mu+\epsilon \psi_1).$

\subsection{Non commutative case}

Assume that $(A,\mu)$ is an associative algebra without commutativity conditions. Let $\mu_t=\mu +t \varphi_1+t^2 \varphi_2+ \cdots$ be a weakly associative formal deformation of $\mu$. As above we put
  $$
  \left\{
  \begin{array}{l}
  \{X,Y\}=\mu(X,Y)-\mu(Y,X), \ X\bullet Y=\mu(X,Y)+\mu(Y,X), \\
   \{X,Y\}_{\mu_t}=\mu_t (X,Y)-\mu_t(Y,X), \ X\bullet_{\mu_t} Y=\mu_t(X,Y)+\mu_t(Y,X).
   \end{array}
   \right.$$
  Since $(A,\mu)$ and $(A[[t]],\mu_t)$ are weakly associative, the brackets $\{X,Y\}$ and $\{X,Y\}_{\mu_t}$ are Lie brackets and $(A,\bullet,\{.,.\})$, $(A[[t]],\bullet_{\mu_t},\{.,.\}_{\mu_t})$ are nonassociative Poisson algebras.  Since $\mu$ is associative, the algebra $(A,\bullet)$ is a commutative Jordan algebra that is
  $$(X \bullet Y)\bullet X^2=X\bullet(Y\bullet X^2)$$
  and $(A[[t]],\bullet_t)$ is a formal deformation of this Jordan algebra.  We have
  \begin{enumerate}[(a)]
  \item order $0$: $\{X,Y \bullet Z\}-\{X,Y\}\bullet Z-Y \bullet\{X,Z\}=0$ and we find again that $(A,\bullet,\{.,.\})$ is a Jordan-Poisson algebra that is the commutative multiplication $\bullet$  is a Jordan multiplication. This shows that for any associative algebra $(A,\mu)$, we can associate a Jordan-Poisson algebra $(A,\bullet,\{.,.\})$ corresponding to the polarized operations. 
    \item order $1$: $B_1(X,Y \bullet Z)+\{X,\rho_1(Y,Z)\}-\rho_1(\{X,Y\},Z)-B_1(X,Y)\bullet Z-Y \bullet B_1(X,Z)-\rho_1(Y,\{X,Z\})=0.$
\end{enumerate}
This implies that we can formulate a notion of non commutative Poisson algebra as follows:
 \begin{definition}
 Let $(A,\mu)$ a non commutative associative algebra and let $(A,\bullet,\{.,.\})$ its polarized Jordan-Poisson algebra. A non commutative Poisson structure on $A$ is given by $(A,\rho_1,  B_1)$ where $B_1$ is a Lie bracket modulo  $B^2_{CE}(A,\{.,.\})$ and $\rho_1$ a Jordan product modulo $B^2_J(A,.\bullet .)$ such that the non commutative Leibniz identity holds:
 \begin{equation}
\label{ncp}
B_1(X,Y \bullet Z)-B_1(X,Y)\bullet Z-Y \bullet B_1(X,Z)+\{X,\rho_1(Y,Z)\}-\rho_1(\{X,Y\},Z)-\rho_1(Y,\{X,Z\})=0.
\end{equation}
 \end{definition}
 In particular, any formal deformation $(A[[t]],\mu_t)$ of the associative algebra $(A,\mu)$ determines a non commutative Poisson structure on $A$ and this deformation could be called a deformation quantization of the non commutative associative algebra $(A,\mu)$.
 
\subsection{Deformations of the Lie bracket in a weakly associative algebra}
Let $(A,\mu)$ a weakly associative algebra and let us consider a formal deformation
$\mu_t=\mu+t \varphi_1+\cdots$ of $\mu$. The following notations are those of the previous section and for simplify we write $\mu(X,Y)=XY$.

We assume in this section that the deformation $\mu_t$ remains invariant the commutative product $X \bullet Y$, that is  $X \bullet_t Y=X \bullet Y$ for any $X,Y$. In this case the bilinear maps $\varphi_i$ are skew-symmetric. The $t$-expansion of the relation $\mathcal{WA}_{\mu_t}=0$ gives:
\begin{enumerate}[(a)]
  \item Order $0$: $\mathcal{WA}(\mu)=0$,
  \item Order $1$: $$
\begin{array}{l}
\varphi_1(X,Y)Z+\varphi_1(XY,Z)-\varphi_1(X,YZ)-X\varphi_1(Y,Z)+\varphi_1(Y,Z)X+\varphi_1(YZ,X)\\
-\varphi_1(Y,ZX)-Y\varphi_1(Z,X) -\varphi_1(Y,X)Z-\varphi_1(YX,Z)+\varphi_1(Y,XZ)+Y\varphi_1(X,Z)=0
\end{array}$$
\end{enumerate}
Since $\varphi_1$ is skew-symmetric
$$
\begin{array}{l}
\varphi_1(\{X,Y\},Z)+\varphi_1(\{Z,X\},Y)+\{\varphi_1(Y,Z),X\}
+2\varphi_1(X,Y)Z-2Y\varphi_1(Z,X)\\
-2\varphi_1(X,YZ)\}
 =0
\end{array}$$
which is equivalent to
$$\begin{array}{l}
\delta^2_L(\varphi_1)(X,Y,Z)+\varphi_1(Y \bullet Z,X)- Y \bullet \varphi_1(Z,X)-Z \bullet \varphi_1(Y,X) =0.
 \end{array}$$
From (\ref{ncp}), since $\rho_1=0$ and $B_1=2\varphi_1$, we have
$$\varphi_1(Y \bullet Z,X)- Y \bullet \varphi_1(Z,X)-Z \bullet \varphi_1(Y,X) =0$$
that implies
$$\delta^2_L(\varphi_1)=0.$$

\begin{proposition}
Let $(A,\bullet,\{.,.\})$ be the nonassociative Poisson algebra which is the polarized of a weakly associative algebra $(A,\mu)$. 
Let $\mu_t=\mu+t \varphi_1+\cdots $ be a formal  deformation of $\mu$ such that $\bullet_t=\bullet$. Then  $\varphi_1$ is a $2$-cocycle of the Lichnerowicz cohomology of the Poisson algebra  $(A,\bullet,\{.,.\})$.
\end{proposition}

\section{Homology of Free weakly associative algebras}

In the first section, we have described the free weakly associative algebra with one generators  $\mathcal{F}_{WA}(X)$. we have seen also that any weakly associative algebra is flexible, that is the associator satisfies the identity
$$\mathcal{A}(X,Y,X)=0.$$
We deduce
\begin{proposition}
The free weakly associative algebra with one generators  $\mathcal{F}_{WA}(X)$ coincides with the free flexible algebra with one generator.
\end{proposition}
This is not true for free weakly associative algebras with $n$ generators, $n >1$. Since these algebras are flexible, there are isomorphic to a quotient of the correspondent free flexible algebras. Recall also that the operad of flexible algebras is not Koszul. In particular, free flexible algebras are not Koszul algebras. This remark conduces to study free weakly associative algebras with the homological point of view.

Concerning the associative case, the free algebra of one generator $\mathcal{F}_{A}(X\}$ is isomorphic to the abelian algebra $\K[X]$ of polynomials of one indeterminate with coefficients in $\K$. Its homology is trivial except in degree $0$. The Hochschild homology for associative algebra corresponds to the complex $(C_n(A,A)=A \otimes A^{\otimes^n}, b_n)$ where
$$\begin{array}{ll}
b_n(m,a_1,\cdots,a_n)= &(ma_1,a_2,\cdots,a_n)+\sum (-1)^i(m,a_1,\cdots a_ia_{i+1},\cdots,a_n)\\
&+(-1)^n(a_nm,a_1,\cdots,a_{n-1}).
\end{array}$$
Let us note that $b_1b_2(a,b,c)=\mathcal{A}(a,b,c)+\mathcal{A}(b,c,a)+\mathcal{A}(c,a,b)$ and the associativity implies that $b_1b_2=0.$ But we have also
$$b_1b_2(a,b,c)=\mathcal{WA}(a,b,c)+\mathcal{WA}(b,a,c)+\mathcal{WA}(c,a,b).$$
Then the weak associativity implies also $b_1b_2=0.$ 

To define an homology for weakly associative algebra, we consider also the complex $$(C_n(A,A)=A \otimes A^{\otimes^n}, b_n^{WA})$$
with $$b_1^{WA}=b_1$$
and $$b_2^{WA}(a_1,a_2,a_3)=b_2(a_1,a_2,a_3)+b_2(a_2,a_3,a_1)-b_2(a_2,a_1,a_3).$$
It is clear that $b_1^{WA}b_2^{WA}=0.$ But, since the free algebra $\mathcal{F}_{WA}(X)$ is commutative, we deduce that $b_2$ and $b_2^{WA}$ coincides on this algebra.  
\begin{proposition}The first homological space $H_0(\mathcal{F}_{WA}(X))$ is isomorphic to $\mathcal{F}_{WA}(X)$.
\end{proposition}
\pf In fact, this algebra is commutative. 

Let us compute 
$H_1(\mathcal{F}_{WA}(X))=\ker b_1/\im \ b_2=C_1(A,A)/\im \ b_2$ with $A=\mathcal{F}_{WA}(X).$  But $C_1(A,A)=A \otimes A= \oplus_{k \geq 1}C_1^k(A,A)$, this grading being defined by the degree of the chains. We have in particular
$$\dim C_1^{2k+1}(A,A)=4d_{2k+1}, \ \ \ \dim C_1^{2k}(A,A)=4d_{2k}-d_k.$$
Likewise $C_2(A,A)=\oplus_{k \geq 0} C_2^k(A,A)$ and if $b_2^k$ denotes the restriction of $b_2$ to $C_2^k(A,A)$ then 
$b_2^k(C_2^k(A,A)) \subset C_1^k(A,A)$. 
Let $H_1^k(A,A)=C_1^k(A,A) / \im b_2^k$ be the homogeneous component correspond to the grading of $C_1(A,A)$.  We have
\begin{enumerate}
\item $\dim H_1^0(A,A)=0.$ In fact $b_2(1,1,1)=(1,1)$ and $\im b_2^0=C_2^0(A,A)=\ker b_1^0.$
  \item $\dim H_1^1(A,A)=1$. In fact $\im b_2^1= \K\{(X,1)\}$ and $\ker b_1^1=C_1^1(A,A)=\K\{(X,1),(1,X)\}.$
  \item Likewise, $\dim H_1^2(A,A)=\dim H_1^3(A,A)=\dim H_1^4(A,A)=1$
  \item $\dim H_1^5(A,A)=2,$
  \item $\dim H_1^6(A,A)=5$
\end{enumerate}
For $k=6m$, we have
$$\begin{array}{ll}
\dim H_1^{6m}(A,A)=&2d_kd_0^2+2d_{k-2}d_1^2+3d_1(d_2d_{k-3}+\cdots+d_{3m-1}d_{3m})\\
&+2d_{k-4}d_2^2+3d_2(d_3d_{k-5}+\cdots+d_{3m-2}d_{3m})+d_2d_{m-1}(d_{m-1}+1)+d_2d_{m-1}^2\\
&+ \cdots \\
&+2d_{m-1}^2d_{m+2}+3d_{m-1}d_md_{m+1}+d_m^2(d_m+1).
\end{array}
$$

 To define other spaces of homology, we shall remark that $\mathcal{F}_{WA}(X)$ is also the free flexible algebra with one generators. To begin, we remark that
 $$b_2b_3(a_1,a_2,a_3,a_4)=b_2((a_1a_2,a_3,a_4)-(a_1,a_2a_3,a_4)+(a_1,a_2,a_3a_4)-(a_4a_1,a_2,a_3)$$
 implying
 $$\begin{array}{ll}
 b_2b_3(a_1,a_2,a_3,a_4)= & (\mathcal{A}(a_1,a_2,a_3),a_4)-(\mathcal{A}(a_4,a_1,a_2),a_3)+(\mathcal{A}(a_3,a_4,a_1),a_2)\\&+(a_1,\mathcal{A}(a_2,a_3,a_4)).
 \end{array}$$
 Let us put
 $$b_3^{WA}(a_1,a_2,a_3,a_4)=b_3(a_1,a_2,a_3,a_4)+b_3(a_1,a_4,a_3,a_2).$$
 Then
 $$b_2b_3^{WA}=b_2b_3\circ\Phi^A_{Id+\tau_{13}}.$$
 But the flexibility implies that $\mathcal{A}\circ \Phi^A_{Id+\tau_{13}}=0$ giving
 $$b_2b_3^{WA}=b_2^{WA}b_3^{WA}=0.$$
 
 Let us compute $H_2(\mathcal{F}_{WA}(X))=\ker b_2/\im \ b_3^{WA}$ with $A=\mathcal{F}_{WA}(X).$ We consider the natural grading of 
$C_3(A,A)=\oplus_{k \geq 0}C_3^k(A,A)$ and the restriction $(b_3^{WA})^k$ to $C^3_k(A,A)$. We have
\begin{enumerate}
  \item $ \im (b_3^{WA})^0= \ker b_2^0=\{0,0,0\}$
  \item $\dim H_2^1(A,A)=1$. In fact, for any $U \in\mathcal{F}_{WA}(X)$, we have
  $$b_3^{WA}(U,1,1,1)=b_3^{WA}(1,1,U,1)=0,$$
  $$b_3^{WA}(1,U,1,1)=b_3^{WA}(1,1,1,U)=(1,1,U)-(1,U,1).$$
  This implies that $\dim (b_3^{WA})^1= 1.$ But $\ker(b_2)^1=\K\{(X,1,1)-(1,X,1),(1,1,X)-(1,X,1)\}$ and $\dim \ker b_2=2.$ Then
  $$\dim H_2(\mathcal{F}_{WA}(X))^1=1.$$
  \item $\dim H_2^2(A,A)=2$. In fact $\im (b_3^{WA})^2=$ is generated by $(1,1,X^2)-(1,X^2,1),(X,X,1)-(X,1,X)$. But
  $$\dim \ker(b_2^2)=\dim \mathcal{C}_2^2(A)-\dim \Im b_2^2=6-2=4.$$
  Then
  $$\dim H_2(\mathcal{F}_{WA}(X))^2=2.$$
  From these results, we deduce
  \begin{proposition}
  The free weakly associative algebra with one generator is not a Koszul algebra.
  \end{proposition}
  \end{enumerate} 
	
	\medskip
	
	\noindent{\bf Remerciements} Je tiens tout particulièrement à remercier Kévin Morand pour toutes les discussions que nous avons eues par courrier, email, skype, bref tous les moyens modernes de communiquer sans se déplacer, et pour m'avoir orienté sur ce sujet.

\end{document}